\documentclass{elsart3-1}

\usepackage{amssymb}

\usepackage[english,francais]{babel}

\newcommand{\R}{\mathbb R}

\newcommand{\HP}{\mathbb H}

\newcommand{\N}{\mathbb N}

\newcommand{\Q}{\mathbb Q}
\newcommand{\Z}{\mathbb Z}

\newcommand{\SI}{\mathbb S}

\newcommand{\GGamma}{{\rm I}\!\Gamma}

\newcommand{\gzg}[1]{[\! [ #1 ]\! ]}

\newcommand{\bast}{{}^{\ast}}
\newcommand{\bstar}{{}^{\bullet}}
\newcommand{\bcirc}{{}^{\circ}}

\newcommand{\fgrm}{\gzg{\pi}_{1}}

\newtheorem{theorem}{Theorem}[section]

\newtheorem{e-proposition}[theorem]{Proposition}

\newtheorem{conjecture}[theorem]{Conjecture}
\newtheorem{e-definition}[theorem]{Definition\rm}

\newtheorem{question}{\bf Question\/}
\newtheorem{ffact}{\bf Fact\/}

\renewcommand{\theequation}{\arabic{equation}}
\def\og{\leavevmode\raise.3ex\hbox{$\scriptscriptstyle\langle\!\langle$~}}
\def\fg{\leavevmode\raise.3ex\hbox{~$\!\scriptscriptstyle\,\rangle\!\rangle$}}

\journal{the Acad\'emie des sciences}
\begin{document}
\centerline{}
\begin{frontmatter}


\selectlanguage{english}
\title{The External Fundamental Group of an Algebraic Number Field}


\selectlanguage{english}
\author{T.M. Gendron}
\ead{tim@matcuer.unam.mx}

\address{Instituto de Matem\'{a}ticas \\
Universidad Nacional Aut\'{o}noma de M\'{e}xico \\
Av. Universidad s/n, Lomas de Chamilpa \\
Cuernavaca CP , Morelos \\ M\'{e}xico}


\medskip

\begin{abstract}
\selectlanguage{english}
We associate to every algebraic number field $K/\Q$ a hyperbolic surface lamination and an {\it external fundamental group} $\bcirc \GGamma_{K}$: a generalization of the fundamental germ construction of  \cite{G1}, \cite{G2} that necessarily contains external (not first order definable) elements. The external fundamental group $\bcirc \GGamma_{\Q}$ is an extension of the absolute Galois group $\hat{\GGamma}_{\Q}$, that conjecturally contains a subgroup whose abelianization is isomorphic
to the id\`{e}le class group.

\vskip 0.5\baselineskip

\selectlanguage{francais}
\noindent{\bf R\'esum\'e} \vskip 0.5\baselineskip \noindent
{\bf Le Groupe Fondamental Externe d'un Corps de Nombres Algebriques}
On associe \`{a} chaque corps de nombres algebriques $K/\Q$ une lamination par surfaces hyperboliques et un {\it groupe fondamental externe} $\bcirc \GGamma_{K}$:
une g\'{e}n\'{e}ralisation de la construction du germe fondamental de \cite{G1},
\cite{G2}, qui contient n\'{e}cessairement des \'{e}l\'{e}ments externes (non definissables au premier ordre). Le groupe fondamental externe $\bcirc \GGamma_{\Q}$ est une extension du groupe de Galois absolut $\hat{\GGamma}_{\Q}$, qui contient conjecturalmente un sous groupe avec une abelianisation isomorphe
au groupe des classes des id\`{e}les.

\end{abstract}
\end{frontmatter}

\selectlanguage{francais}

\selectlanguage{english}
\section{Introduction}
\label{}

The search for a geometrization of an algebraic number field $K/\Q$ has been one of the longstanding ambitions
of algebraic number theory:  indeed, it could be said that the specter of such a 
geometrization
haunts some of its most celebrated enterprises {\it viz.}\
the Riemann hypothesis, nonabelian class field theory, 
Grothendieck-Teichm\"{u}ller theory.   One phenomenon which could achieve structural clarity via geometrization is the isomorphism of class field theory $C_{\Q}\cong \R^{\times}_{+}\times\hat{\Z}^{\times}$,
where $C_{\Q}$ is the idele class group of $\Q$.  Since $\hat{\Z}^{\times}\cong
\hat{\GGamma}^{\rm ab}_{\Q}$, where $\hat{\GGamma}_{\Q}={\rm Gal}(\bar{\Q}/\Q )$, it has been suggested by a number of authors \cite{We}, \cite{Ca}, \cite{Co} that the
factor $ \R^{\times}_{+}$ ought to have also a Galois interpretation.  Formally, one seeks an extension $\bar{\GGamma}_{\Q}\rightarrow \hat{\GGamma}_{\Q}$ in which
$\bar{\GGamma}_{\Q}$ has arithmetic meaning (a ``cosmic Galois group" \cite{Ca}), and  
for which $\bar{\GGamma}^{\rm ab}_{\Q}\cong C_{\Q}$.  
In this paper, we shall construct a candidate for $\bar{\GGamma}_{\Q}$, defined as the external fundamental group of a geometrization of $\Q$ by a hyperbolic surface lamination.  

{\it Acknowledgements :} I would like to thank P. Lochak, who made many useful comments
during the initial phases of this work.  In addition, I would like to thank the International Centre for Theoretical Physics, which supported
a visit for a month during the summer of 2007, when much of this work was carried out.  This work was supported by the author's PAPIIT grant IN103708  and CONACyT grant 58537.

\section{Internal Fundamental Group}\label{secfundgerm}

Let $M$ be a compact $n $-manifold, $p:\tilde{M}\rightarrow M$ a universal cover
and write $\pi=\pi_1(M)$.  Fix an ultrafilter $\mathfrak{U}$
on $\N$ all of whose elements are of infinite cardinality.  Denote by $\bast\pi$ the ultraproduct
of $\pi$ with respect to $\mathfrak{U}$.  Note that there is a monomorphism 
 $c:\pi\hookrightarrow\bast\pi$ given by the constant sequences, and we identify $\pi$
 with its image.  The ultraproduct $\bast\pi$ is an example of a nonstandard model of $\pi$ \cite{Ro}.

Suppose that $M$ is riemannian, and equip $\tilde{M}$ with the pull-back metric so that
$\pi$ acts by isometries on $\tilde{M}$.   Let $\bstar\tilde{M}$ be the quotient of $\bast\tilde{M}$ 
(= the ultraproduct of $\tilde{M}$) obtained
by identifying sequence classes that are asymptotic.  There is a canonical surjective
 map $\bstar\tilde{M}\rightarrow M$ which associates to each class $\bstar\tilde{x}$ the limit of $p(\tilde{x}_{i})$
 where $\{ \tilde{x}_{i}\}\in \bstar\tilde{x}$ is any representative sequence for which  $p(\tilde{x}_{i})$ converges.  Note that $\bast\pi$ acts on the left on $\bstar \tilde{M}$ and $\bast\pi\backslash \bstar\tilde{M}\approx M$.

We may view $\bstar\tilde{M}$ as a lamination with discrete
transversals: the leaf containing $\bstar \tilde{x}\in\bstar\tilde{M}$ consists of those sequence classes of bounded distance from $\bstar \tilde{x}$, 
itself a riemannian manifold.  In fact, $\bstar\tilde{M}$ may be identified with the suspension of the inclusion $c$, {\it i.e.}\
$ \left( \tilde{M}\times\bast\pi\right)/ \pi $,
where $ (\tilde{x},\bast \alpha ) \cdot \gamma = (\gamma\cdot\tilde{x}, (\bast\alpha )\gamma^{-1})$
for all $\gamma\in\pi$.
 In the suspension description,
the action of $\bast\pi$ is induced by $ (\tilde{x},\bast \alpha ) \mapsto  (\tilde{x},\bast\gamma \bast\alpha )$, $\bast\gamma\in\bast\pi$, and so can be seen to be by leaf-wise isometries.
This discussion  applies to any group extension $\pi \subset G$ (particularly, when
$G$ = a nonstandard model of $\pi$),  
the appropriate universal covering space being the suspension of the inclusion 
$\pi\hookrightarrow G$.

We now indicate how $\bast\pi$ codifies laminated
coverings of $M$.  For simplicity, we shall restrict ourselves to suspensions over $M$. 
 Let $G$ be a compact topological group
and let $\rho :\pi\rightarrow G$ be a representation.  The suspension of $\rho$, denoted $M(\rho )$, is a principal 
$G$-bundle as well as a lamination over $M$, minimal if and only if $\rho$ has dense image, with simply connected leaves
if and only if ${\rm Ker}(\rho )=1$.  Three examples:

\begin{enumerate}
\item[a.]  If $G=1$ then $M(\rho )\approx M$.
\item[b.]  Let $G=\hat{\pi}$ =  the profinite completion of $\pi$, $\rho$ the canonical
map. Then $M(\rho )$ $\approx$ $\hat{M}$ = the algebraic universal cover of $M$, a $\hat{\pi}$-principal
bundle over $M$ {\it e.g.}\ $\hat{\pi}\backslash \hat{M}\approx M$. It is classical that
$\hat{M}$ and $\hat{\pi}$ are the appropriate notions of universal cover and fundamental group  
for $M$ within the \'{e}tale category. 
\item[c.] Let $M=G=\SI^{1}=\R/\Z$, and for $r\in\R-\Q$, define $\rho$ by $\rho (n )=
\overline{nr}$ = the image of $nr$ in $\SI^{1}$.  Then $M(\rho ) = \mathcal{F}_{r}$ = 
the irrational foliation of the 2-torus 
by lines of slope $r$.
\end{enumerate}

An analogue of the fundamental group for $M(\rho )$ is given by the {\it fundamental germ}
$\gzg{\pi}=\fgrm M(\rho )$, \cite{G1}, \cite{G2}.  In the case when the suspension $M(\rho )$ is minimal, it has the following description.    Since $\rho$ has dense image,
the ``standard part'' map $ {\rm std}(\rho ): \bast\pi \rightarrow G  $,
defined by taking a sequence class to the unique limit in $G$ of its image by $\rho$, is onto.  We define 
$ \gzg{\pi}  := {\rm Ker}( {\rm std}(\rho ))$
and refer to $1\rightarrow \gzg{\pi}\rightarrow\bast\pi\rightarrow G\rightarrow1$
as the {\it standardization exact sequence}.  For the three examples above we have: 

\begin{enumerate}
\item[a.]  $\fgrm M=\bast\pi$.
\item[b.] $\fgrm\hat{M}= \bigcap \bast H$ where $H<\pi$ runs through
the subgroups
of finite index.  This is a non-trivial subgroup of $\bast\pi$ even when $\pi$ is residually
finite {\it i.e} when $\bigcap H$ is trivial (for example, when $M$ is a compact surface). 
\item[c.]  We say that a sequence class 
$\bast\epsilon\in\bast\R$
is an infinitesimal if it contains a sequence converging to 0.  
Then we may identify $\fgrm\mathcal{F}_{r}$ with the subgroup of $\bast n\in\bast\Z$ for which $r\bast n +\bast m$ is an infinitesimal for some $\bast m\in \bast\Z$: in other words, $\fgrm\mathcal{F}_{r}$ is the group
of diophantine approximations of $r$.  
 \end{enumerate}
 
 We now discuss covering space theory. 
 Let $M(\rho )$ be as above, assumed for simplicity to be minimal with simply connected
 leaves.  Assume also that $M$ has been equipped with a
riemannian metric, so that $M(\rho )$ has a leaf-wise riemannian metric.  
There is a canonical map $\tilde{M}\rightarrow M(\rho )$, induced by $\tilde{M}\times \{ 1\}\hookrightarrow \tilde{M}\times G$.  The image of this map is a leaf $L_{0}$ called the canonical
leaf.
There is a surjective map $\bstar \tilde{M}\rightarrow M(\rho )$ -- assigning to a sequence
class $\bstar\tilde{x}$ the limit of its image via $\tilde{M}\rightarrow M(\rho )$ -- which is a local isometry along the leaves.  Any continuous self-map of $M(\rho )$ preserving $L_{0}$ 
lifts uniquely to a self-map of $\bstar \tilde{M}$.
The natural action of $\gzg{\pi}$ on  $\bstar\tilde{M}$ has quotient 
$\gzg{\pi}\backslash\bstar\tilde{M}$
which is in canonical bijection with $M(\rho )$.
For example, when $M=\Gamma\backslash\HP^{2}$ is a closed hyperbolic surface, we may identify 
$\fgrm M(\rho )$ with a ``Fuchsian germ'' $\gzg{\Gamma}<{\rm PSL}(2,\bast\R )$ and $\gzg{\Gamma}\backslash\bstar\HP^{2}$ is in bijection with $M(\rho )$.  

It is possible to endow $\bstar\tilde{M}$ with a nontrivial transverse topology in such a way
that $\gzg{\pi}$ acts by homeomorphisms and so that the quotient $\gzg{\pi}\backslash \bstar\tilde{M}$ is homeomorphic to $M (\rho )$. 
To do this, we choose a set-theoretic section $s:G\rightarrow \bast \pi$ of ${\rm std}(\rho )$, so that $s(\rho (\gamma ))=\gamma$ for all $\gamma\in\pi$,
and for which $s(G)$ is a right $\pi$-set.  Then if we give $\bast\pi$ the topology:
(topology of $G$) $\times$ (discrete), this gives a topology
on $\tilde{M}\times \bast\pi$ invariant by the action of $\pi$, hence inducing
a topology on $\bstar\tilde{M}$.  The left multiplication action by elements of $\gzg{\pi}$ permutes the
``cosets'' $\gzg{x}s(G)$, $\gzg{x}\in\gzg{\pi}$, hence $\gzg{\pi}$ acts by homeomorphisms,
and with the quotient topology, the bijection between $\gzg{\pi}\backslash \bstar\tilde{M}$ 
and $M(\rho )$
is a homeomorphism.
(N.B. We may even choose the section $s$ in order that any leaf of $\bstar\tilde{M}$ intersects a
given $s(G)$-transversal
no more than once: so that $\bstar\tilde{M}$ is a  lamination with no non-trivial holonomy.)

\section{External Fundamental Group}

Let $F$ be the free group on two generators, $\hat{F}$ its profinite completion
and consider the standardization sequence $1\rightarrow \gzg{F}\rightarrow \bast F\rightarrow \hat{F}\rightarrow 1$.  Neither $\bast F$ nor $\hat{F}$ are free groups in the discrete (combinatorial)
sense.  Let $\hat{{\bf F}}$ be the free group generated by $\hat{F}$ (viewed as a set), which has cardinality of the continuum.   By universality, there is a canonical epimorphism 
$\hat{p}:\hat{\bf F}\rightarrow \hat{F}$.  If $\sigma :\hat{F}\hookrightarrow \bast F$ is a set-theoretic section
of the standardization sequence whose image contains a generating set of $\bast F$, then the induced map $\bast p: \hat{\bf F}\rightarrow \bast F$
is an epimorphism, and $\hat{p}={\rm std}\circ \bast p$ (by the uniqueness part of universality).  If $\hat{K}$, $\bast K$ are the kernels
of $\hat{p}, \bast p$ then $\bast K<\hat{K}$.

Denote by ${\rm Aut}(\hat{F})$ the group of bicontinuous automorphisms of $\hat{F}$,
and by $\bcirc {\rm Aut}(F)$ the subgroup of ${\rm Aut}(\bast F )$ of automorphisms which
induce elements of ${\rm Aut}(\hat{F})$  {\it i.e.}\ automorphisms which stabilize $\gzg{F}$ 
and induce bicontinuous automorphisms of $\hat{F}$.  Note
that ${\rm Aut}(F)$ as well as $\bast {\rm Aut}(F)$ include canonically in $\bcirc {\rm Aut}( F )$.
Indeed, if $\bast A\in \bast {\rm Aut}(F)$ and $\bast x\in\gzg{F}$, then $\bast A (\bast x )$
is represented by a sequence $\{ A_{i}(x_{i})\}$, and $A_{i}(x_{i})$ is in a subgroup of
index $N_{i}\rightarrow\infty$ if and only if $x_{i}$ is.  

\begin{theorem}\label{section} The canonical homomorphism
$\bcirc {\rm Aut}(F)\rightarrow {\rm Aut}(\hat{F})$ is surjective.  \end{theorem}

The Theorem is proved as follows: note first that any element $\alpha\in  {\rm Aut}(\hat{F})$ defines a bijection of  the generating set of $\hat{\bf F}$, hence an automorphism {\boldmath$\alpha$} of the latter.
As such,  {\boldmath$\alpha$} necessarily stabilizes $\hat{K}$: we may arrange that it also
stabilizes $\bast K$ by
composing, if necessary, with a suitable automorphism covering the identity of $\hat{F}$.
The result descends to an automorphism $\bcirc\alpha$ of $\bast F$. The association 
$\alpha\mapsto \bcirc\alpha$ evidently defines a (set-theoretic) section.

Denote by $\bcirc {\rm Inn}(F)$ those elements of $\bcirc {\rm Aut}(F)$ which
map to inner automorphisms of $\hat{F}$.  (N.B. $\bast F$, acting
innerly, is a subgroup of $\bcirc {\rm Inn}( F)$.)  If we denote by $\bcirc {\rm Out}( F)$
the quotient of $\bcirc {\rm Aut}(F)$ by  $\bcirc {\rm Inn}(F)$, we obtain an exact sequence
$1\rightarrow \gzg{\GGamma}\rightarrow \bcirc {\rm Out}( F) \rightarrow {\rm Out}(\hat{F})\rightarrow 1 $. 

It is important to note that $\bcirc {\rm Out}(F)$ contains as a {\it proper} subgroup
the ultraproduct $\bast {\rm Out}(F)\cong \bast {\rm GL}(2,\Z )\cong {\rm GL}(2,\bast \Z )$.  The latter is called
the group of {\it internal} outer automorphisms of $\bast F$, and elements of $\bcirc {\rm Out}(F)$
which are not internal are called {\it external}.  That we cannot replace
$\bcirc {\rm Out}( F )$ by $\bast {\rm Out}(F)$ is borne out by the following
\begin{ffact}  Although $F$ is dense in $\hat{F}$, ${\rm Out}(F)$ is not dense
in ${\rm Out}(\hat{F})$, hence ${\rm Out}(\hat{F})$ is not the profinite completion of
${\rm Out}(F)\cong {\rm GL}(2,\Z )$.  Thus,  $\bast {\rm Out}(F)$
does not map epimorphically onto ${\rm Out}(\hat{F})$.
\end{ffact}

Recall that the theory of a group $G$ is the collection ${\rm Th}(G)$ of all first order sentences
which are true in $G$.  We say $G'$ is a nonstandard model of $G$ if 
${\rm Th}(G')={\rm Th}(G)$ but $G'\not\cong G$.  For example, the ultrapower $\bast G$ is a
nonstandard model of $G$.

\begin{question}  Is  $\bcirc {\rm Out}( F )$ a nonstandard model of ${\rm Out}(F)$?
\end{question}

In what follows $K/\Q$ is an arbitrary algebraic number field and $\hat{\GGamma}_{K}$
is its absolute Galois group.
Recall the Belyi monomorphism
$\beta :\hat{\GGamma}_{K}\subset \hat{\GGamma}_{\Q}\hookrightarrow {\rm Out}(\hat{F})$.
We will not distinguish between $\hat{\GGamma}_{K}$ and its image in $ {\rm Out}(\hat{F})$.
Let ${\rm SL}(2,\Z )\cong {\rm Out}_{+}(F)\hookrightarrow {\rm Out}(\hat{F})$ be the canonical
inclusion.  Define $\hat{\Sigma}=\hat{\Sigma}_{\bar{\Q}}$ as the suspension
$  (\HP^{2}\times {\rm Out}(\hat{F}) )/{\rm SL}(2,\Z )$,
where the action of $A\in {\rm SL}(2,\Z )$ is defined $A(z, f) = (Az, fA^{-1})$.  Then
$\hat{\Sigma}$ is a non-minimal solenoid by hyperbolic surface orbifolds that covers the modular orbifold ${\rm SL}(2,\Z )\backslash \HP^{2}$.  The action of $\hat{\GGamma}_{K}$ on the product
$\HP^{2}\times {\rm Out}(\hat{F})$, $\hat{\sigma} (z,f)= (z,\hat{\sigma}f)$, descends to
an action on $\hat{\Sigma}$ by leaf-wise isometries.  Since $\hat{\GGamma}_{K}$
is a closed subgroup of ${\rm Out}(\hat{F})$, the quotient
$\hat{\Sigma}_{K} = \hat{\GGamma}_{K}\backslash \hat{\Sigma}$
is also a lamination by hyperbolic surface orbifolds.  By construction, the association
$K\mapsto \hat{\Sigma}_{K}$ is Galois natural.

 Denote by $ \bcirc{\GGamma}_{K}$ the pre-image of $\hat{\GGamma}_{K}$ in
 $\bcirc {\rm Out} (F )$ so that $\gzg{\GGamma}$ 
is the kernel of the standardization
 epimorpism $\bcirc{\GGamma}_{K}\rightarrow \hat{\GGamma}_{K}$.  
 We
 have $\gzg{\GGamma}=\bigcap\bcirc\GGamma_{K}$.
Recall that there is a canonical inclusion 
${\rm SL}(2,\Z )\cong {\rm Out}_{+}(F)\hookrightarrow \bcirc {\rm Out}(F )$. 
By suspending this inclusion with respect to the action of
${\rm SL}(2,\Z )$ on $\HP^{2}$, we obtain a trivial lamination
which we denote $\bcirc\HP^{2}$.  We note that the quotient of $\bcirc\HP^{2}$ 
by the left action of $ \bcirc {\rm Out}(F )$
is isometric to the modular orbifold.

We topologize $\bcirc {\rm Out}(F)$ by choosing a set-theoretic section 
of $\bcirc {\rm Out}(F)\rightarrow {\rm Out}(\hat{F})$ whose image
is a right ${\rm SL}(2,\Z )$-set and which maps ${\rm SL}(2,\Z )$ to itself
(as we did at the end of the last section).  This induces a topology
on $\bcirc\HP^{2}$ making it a solenoid by hyperbolic surface orbifolds, with respect to which
the action by $\bcirc {\rm Out}(F)$ is by homeomorphisms which are isometries along
the leaves.  The quotient by $\gzg{\GGamma}$ can be identified with $\hat{\Sigma}=\hat{\Sigma}_{\bar{\Q}}$ 
and in addition
$   \hat{\Sigma}_{K} \cong  \bcirc \GGamma_{K} \backslash\bcirc\HP^{2}\cong\hat{\GGamma}_{K}
\backslash \hat{\Sigma}$.
This justifies viewing  $\bcirc\GGamma_{K}$ as a fundamental group, in a way which 
generalizes the internal fundamental group defined in \S 2.   

\begin{conjecture}  There is a subgroup $\bar{\GGamma}_{\Q}<\bcirc\GGamma_{\Q}$ which
is an extension of $\hat{\GGamma}_{\Q}$ and for which $\bar{\GGamma}_{\Q}^{\rm ab}\cong
C_{\Q}$.
\end{conjecture}





\end{document}